# THE PFAFFIAN CLOSURE OF AN O-MINIMAL STRUCTURE

PATRICK SPEISSEGGER

ABSTRACT. Every o-minimal expansion $\widetilde{\mathbb{R}}$ of the real field has an o-minimal expansion $\mathcal{P}(\widetilde{\mathbb{R}})$ in which the solutions to Pfaffian equations with definable $C^1$ coefficients are definable.

## INTRODUCTION

The main result of this paper is the following.

**Theorem.** *Let $\widetilde{\mathbb{R}}$ be an o-minimal expansion of the real field. Then there is an o-minimal expansion $\mathcal{P}(\widetilde{\mathbb{R}})$ of $\widetilde{\mathbb{R}}$ which is closed under solutions to Pfaffian equations in the following strong sense (where "definable" refers to definability in $\mathcal{P}(\widetilde{\mathbb{R}})$). Whenever $U$ is a definable and connected open subset of $\mathbb{R}^n$, $\omega = a_1 dx_1 + \cdots + a_n dx_n$ is a 1-form on $U$ with definable coefficients $a_i : U \longrightarrow \mathbb{R}$ of class $C^1$, and $L \subseteq U$ is a Rolle leaf of $\omega = 0$, then $L$ is also definable.*

The notion of Rolle leaf is due to Moussu and Roche [10], who work in the analytic setting and were inspired by the Khovanskii-Rolle Theorem [7].

As a special case, note that if $U \subseteq \mathbb{R}^n$ is as in the theorem and $f : U \longrightarrow \mathbb{R}$ is a $C^1$ function satisfying
$$\frac{\partial f}{\partial x_i}(x) = F_i(x, f(x)), \quad x \in U, \quad i = 1, \ldots, n,$$
where each $F_i : U \times \mathbb{R} \longrightarrow \mathbb{R}$ is a definable $C^1$ function, then $f$ is also definable. (This follows from the theorem because the graph of $f$ is a Rolle leaf of $\omega = 0$, where $\omega = F_1 dx_1 + \cdots + F_n dx_n - dy$ on $U \times \mathbb{R}$; see Section 1.)

**Corollary.** *Suppose that $I \subseteq \mathbb{R}$ is an open interval, $a \in I$ and $g : I \longrightarrow \mathbb{R}$ is definable and continuous. Then its antiderivative $f : I \longrightarrow \mathbb{R}$ given by $f(x) := \int_a^x g(t) dt$ is also definable.*

The proof of the theorem goes as follows. After defining Rolle leaves in the $C^1$ setting and establishing some basic facts about them in Section 1, it is shown in Section 2 that the Khovanskii theory according to Moussu and Roche [10, 11] carries over to the o-minimal setting with $C^1$ data, and is actually simplified in the process. (The analytic data and semianalytic assumptions in [10] make certain precautions necessary which disappear in the o-minimal context.) Section 3 then adapts the more recent work by Lion and Rolin [8] on $T^\infty$-Pfaffian sets to the o-minimal setting. (This is done here without using the "Cauchy-Crofton formula".) The above theorem is then an easy consequence.

*Date*: August 1997; Preliminary Version.
1991 *Mathematics Subject Classification.* Primary 14P10, 58A17; Secondary, 03C99.
*Key words and phrases.* o-minimal structures, Pfaffian systems.
Research supported by The Fields Institute for Research in Mathematical Sciences.

There are some recent results that point in the direction of the theorem. Wilkie's theorem of the complement [14] implies that the expansion of the real field by the Pfaffian functions in the sense of Khovanskii [5] is o-minimal. These Pfaffian functions are analytic and defined on all of $\mathbb{R}^n$. Wilkie's arguments inspired Lion and Rolin [8] to introduce the notion of a $T^\infty$-Pfaffian set in the semianalytic setting and to prove a corresponding o-minimality result. Karpinski and Macintyre [4] showed that for a given o-minimal expansion $\widetilde{\mathbb{R}}$ of the real field, the expansion of $\widetilde{\mathbb{R}}$ by (total) $C^\infty$ functions which are Pfaffian relative to $\widetilde{\mathbb{R}}$ remains o-minimal.

It remains to be seen whether the extra generality of this paper's main theorem – roughly speaking, $C^1$ leaves versus $C^\infty$ functions – is genuine. This issue may be related to the open problem whether the o-minimal structures of [4, 8, 14] and of this paper are model complete in their natural languages.

**Conventions.** Throughout this paper the letters $k, l, m, n, p, q$ range over $\mathbb{N} = \{0, 1, 2, \dots\}$.

For any set $S$ we denote by $|S|$ the cardinality of $S$.

We let $\langle v_1, \dots, v_k \rangle$ denote the linear span of vectors $v_1, \dots, v_k$ in a vector space $E$.

A **box** in $\mathbb{R}^n$ is a cartesian product of open intervals $(a, b)$ with $a, b \in \mathbb{R} \cup \{-\infty, \infty\}$ and $a < b$. If we only allow $a, b \in \mathbb{Q} \cup \{-\infty, \infty\}$, then we call the resulting box a **rational box**.

Let $A \subseteq \mathbb{R}^m$. We write $\text{cl}(A)$, $\text{int}(A)$, $\text{bd}(A) := \text{cl}(A) \setminus \text{int}(A)$ and $\text{fr}(A) := \text{cl}(A) \setminus A$ for the topological closure, interior, boundary and frontier of $A$, respectively. By "component" we always mean "connected component".

For $0 \leq k \leq m$ and any $x \in \mathbb{R}^k$, we denote by $A_x := \{y \in \mathbb{R}^{m-k} : (x, y) \in A\}$ the fiber of $A$ over $x$.

We equip $\mathbb{R}^m$ with the usual metric $d$, which for $m > 0$ is given by $d(x, y) := \max\{|x_1 - y_1|, \dots, |x_m - y_m|\}$. We also set $d(x, A) := \inf_{y \in A} d(x, y)$. For $\delta > 0$ we let $B(a, \delta) := \{x \in \mathbb{R}^m : d(x, a) < \delta\}$ and $T(A, \delta) := \{x \in \mathbb{R}^m : d(x, A) < \delta\}$.

For a map $\iota : \{1, \dots, k\} \longrightarrow \{1, \dots, n\}$ we let $\Pi_\iota : \mathbb{R}^n \longrightarrow \mathbb{R}^k$ be the map given by $\Pi_\iota(x_1, \dots, x_n) := (x_{\iota(1)}, \dots, x_{\iota(k)})$. (In the case $k = 0$ the map $\Pi_\iota$ sends every $x \in \mathbb{R}^n$ to the single element 0 of $\mathbb{R}^0$.)

All manifolds are assumed to be nonempty $C^1$ submanifolds of some $\mathbb{R}^n$, and are also assumed to be embedded submanifolds, unless specifically referred to as immersed submanifolds.

"Definable" always refers to definability with real parameters.

## 1. Generalities on Rolle Leaves

Let $U \subseteq \mathbb{R}^n$ be open and let $\omega = a_1 dy_1 + \cdots + a_n dy_n$ be a **1-form on $U$ of class $C^1$**, that is, each $a_i : U \longrightarrow \mathbb{R}$ is a $C^1$ function. Put

$$S(\omega) := \{\, y \in U : \omega(y) = 0 \,\} = \bigcap_{i=1}^n a_i^{-1}(0)$$

(so $S(\omega)$ is a closed subset of $U$). We call $\omega$ **nonsingular** if $S(\omega) = \emptyset$. We think of "$\omega = 0$" as defining a hyperplane field on $U \setminus S(\omega)$: it assigns to each $y \in U \setminus S(\omega)$ the hyperplane

$$\ker(\omega(y)) = \{\, x \in \mathbb{R}^n : a_1(y)x_1 + \cdots + a_n(y)x_n = 0 \,\} \subseteq \mathbb{R}^n$$

under the usual identification of $T_y(U \setminus S(\omega)) = T_y(\mathbb{R}^n)$ with $\mathbb{R}^n$.

**1.1. Definition.** An **integral manifold** of the equation $\omega = 0$ is an $(n-1)$-dimensional immersed $C^1$ submanifold $M$ of $U \setminus S(\omega)$ such that $T_y M = \ker(\omega(y))$ for all $y \in M$. (Note: since $M \subseteq U \setminus S(\omega)$ is only immersed, the manifold $M$ does not necessarily have the topology induced by $\mathbb{R}^n$.) A **leaf** of $\omega = 0$ is a maximal connected integral manifold of $\omega = 0$.

A leaf $L$ of $\omega = 0$ is a **Rolle leaf** if $L$ is an embedded (not just immersed) submanifold of $U \setminus S(\omega)$ that is closed in $U \setminus S(\omega)$, such that each $C^1$ curve $\gamma : [0, 1] \longrightarrow U$ with $\gamma(0), \gamma(1) \in L$ is tangent at some point to the hyperplane field defined by $\omega = 0$, that is, there is $t \in [0, 1]$ such that $\omega(\gamma(t))(\gamma'(t)) = 0$.

**1.2. Example.** Let $\omega$ be an **integrable** 1-form, that is, $\omega \wedge d\omega = 0$. Then the leaves of $\omega = 0$ form a partition of $U \setminus S(\omega)$, called the **foliation** defined by $\omega = 0$, see [13]. By the Rolle-Khovanskii Lemma [6], if $L$ is a leaf of $\omega = 0$ such that $U \setminus L$ has exactly two connected components and each component has $L$ as frontier in $U$, then $L$ is in fact a Rolle leaf.

*Remark.* Note that Rolle leaves of $\omega = 0$ are defined here even if $\omega$ is not integrable, in contrast to [10]. This is convenient but does not add to the generality of the main result (see [11] for an explanation in the analytic context; a similar argument works here).

**1.3. Example.** Let $V \subseteq \mathbb{R}^n$ be nonempty, open and connected, let $F_1, \ldots, F_n : V \times \mathbb{R} \longrightarrow \mathbb{R}$ be $C^1$ functions, and let $f : V \longrightarrow \mathbb{R}$ be a $C^1$ function satisfying the **Pfaffian system**

$$\frac{\partial f}{\partial x_i}(x) = F_i(x, f(x)) \quad \text{for } x \in V \text{ and each } i = 1, \ldots, n.$$

Then the graph $\Gamma(f)$ of $f$ is a Rolle leaf in $U := V \times \mathbb{R}$ of the equation

$$\omega := F_1 dx_1 + \cdots + F_n dx_n - dy = 0.$$

To see this, note first that $S(\omega) = \emptyset$, and that $\Gamma(f)$ is closed in $U$ and is clearly an embedded $n$-dimensional $C^1$ submanifold of $U$. Next, the connectedness of $V$ implies that $U \setminus \Gamma(f)$ has exactly two connected components $C_1 := \{ (x, y) \in U : y < f(x_1, \ldots, x_n) \}$ and $C_2 := \{ (x, y) \in U : y > f(x_1, \ldots, x_n) \}$. Now let $\gamma : [0, 1] \longrightarrow U$ be a $C^1$ curve with $\gamma(0), \gamma(1) \in \Gamma(f)$. We may clearly assume that $\omega(\gamma(0))(\gamma'(0)) \neq 0$ and $\omega(\gamma(1))(\gamma'(1)) \neq 0$, and that $\gamma((0, 1))$ is contained in one of the two connected components.

We claim that then $\omega(\gamma(0))(\gamma'(0))$ and $\omega(\gamma(1))(\gamma'(1))$ must have different sign. For if $\omega(\gamma(0))(\gamma'(0)) > 0$, say, then there is an $\epsilon > 0$ such that $\gamma((0, \epsilon)) \subseteq C_2$, and so by the above $\gamma((0, 1)) \subseteq C_2$; but if also $\omega(\gamma(1))(\gamma'(1)) > 0$, then there is a $\delta > 0$ such that $\gamma((\delta, 1)) \subseteq C_1$, so that $\gamma((0, 1)) \subseteq C_1$, a contradiction. A symmetric argument works if both $\omega(\gamma(0))(\gamma'(0))$ and $\omega(\gamma(1))(\gamma'(1))$ are negative, and so the claim is proved.

It follows from the claim and Rolle's Theorem that there exists a $t \in (0, 1)$ such that $\omega(\gamma(t))(\gamma'(t)) = 0$, which together with the above proves that $\Gamma(f)$ is a Rolle leaf.

**1.4. Lemma.** *Let $U \subseteq \mathbb{R}^n$ be open and let $\omega$ be a 1-form of class $C^1$ on $U$. Let $L$ be an integral manifold of $\omega = 0$, and suppose that there is a closed set $K \subseteq \mathbb{R}^n$ such that $L = K \cap (U \setminus S(\omega))$. Let $C \subseteq U \setminus S(\omega)$ be a connected manifold of dimension $\leq n - 1$ such that $T_x C \subseteq \ker(\omega(x))$ for all $x \in C$. Then either $C \cap L = \emptyset$ or $C \subseteq L$.*

*Proof.* Write $\omega = a_1 dx_1 + \cdots + a_n dx_n$ with each $a_i : U \longrightarrow \mathbb{R}$ of class $C^1$. Assume that $C \cap L \neq \emptyset$. Since $C \cap L = C \cap K$, the subset $C \cap L$ of $C$ is closed in $C$. By the connectedness of $C$ it is therefore enough to show that $C \cap L$ is also open in $C$.

Fix $x \in C \cap L$. Permuting coordinates if necessary, we may assume that $a_n(y) \neq 0$ for all $y$ in some open neighbourhood $V$ of $x$. Then the projection $\Pi : \mathbb{R}^n \longrightarrow \mathbb{R}^{n-1}$ on the first $n-1$ coordinates is such that $\Pi\big|(V \cap L)$ is an immersion. Hence, after shrinking $V$ if necessary, $V \cap L$ is the graph of a $C^1$ function $f : W \longrightarrow \mathbb{R}$, where $W := \Pi(V)$. This function $f$ satisfies the following equations:

$$(1.1) \qquad \frac{\partial f}{\partial x_i}(x') = -\frac{a_i}{a_n}\big(x', f(x')\big), \quad \text{for } x' \in W \text{ and } i = 1, \ldots, n-1.$$

Let now $C'$ be the connected component of $C \cap V$ that contains $x$; it clearly suffices to show that $C' \subseteq \Gamma(f)$. Let $y \in C'$ be such that $y \neq x$, and choose $\epsilon > 0$ and a $C^1$ curve $\gamma = (\gamma_1, \ldots, \gamma_n) : (-\epsilon, 1+\epsilon) \longrightarrow C'$ such that $\gamma(0) = x$ and $\gamma(1) = y$. Note that by hypothesis on $C$ we have $\gamma'(t) \subseteq \ker(\omega(\gamma(t)))$ for $t \in (-\epsilon, 1+\epsilon)$, that is,

$$\gamma_n'(t) = -\frac{a_1}{a_n}\big(\gamma(t)\big) \cdot \gamma_1'(t) - \cdots - \frac{a_{n-1}}{a_n}\big(\gamma(t)\big) \cdot \gamma_{n-1}'(t)$$
$$= -\frac{a_1}{a_n}\big(\delta(t), \gamma_n(t)\big) \cdot \gamma_1'(t) - \cdots - \frac{a_{n-1}}{a_n}\big(\delta(t), \gamma_n(t)\big) \cdot \gamma_{n-1}'(t)$$

for all such $t$, where $\delta := (\gamma_1, \ldots, \gamma_{n-1}) : (-\epsilon, 1+\epsilon) \longrightarrow W$. On the other hand, we define $h : (-\epsilon, 1+\epsilon) \longrightarrow \mathbb{R}$ by

$$h(t) := f(\delta(t)).$$

Then by (1.1) we have

$$h'(t) = -\frac{a_1}{a_n}\big(\delta(t), h(t)\big) \cdot \gamma_1'(t) - \cdots - \frac{a_{n-1}}{a_n}\big(\delta(t), h(t)\big) \cdot \gamma_{n-1}'(t)$$

for all $t \in (-\epsilon, 1+\epsilon)$. Therefore $h$ and $\gamma_n$ satisfy the same differential equation, and since $h(0) = f(\delta(0)) = x_n = \gamma_n(0)$, if follows from the uniqueness of solutions to ordinary differential equations that $\gamma_n = h$. In particular $y = \gamma(1) = (\delta(1), h(1)) \in \Gamma(f)$, and since $y \in C'$ was arbitrary, this finishes the proof of the lemma. $\square$

**1.5. Definition.** Let $\Omega = (\omega_i)_{i \in I}$ be a finite family of 1-forms of class $C^1$ on an open set $U \subseteq \mathbb{R}^n$, and let $N$ be a submanifold of $U$. We say that $\Omega$ is **transverse** to $N$ if for each $y \in N$ the family $\big(\omega_i(y)\big|T_y N\big)_{i \in I}$ is a linearly independent family of linear forms on $T_y N$, or equivalently, if

$$\dim\left(T_y N \cap \bigcap_{i \in I} \ker(\omega_i(y))\right) = \dim(N) - |I|$$

for all $y \in N$. (Note that then in particular $|I| \leq \dim(N)$ and $S(\omega_i) \cap N = \emptyset$ for each $\omega_i$.) If $|I| = 1$ and $\omega$ is the unique element of $\Omega$, we sometimes write "$\omega$ is transverse to $N$" in place of "$\Omega$ is transverse to $N$".

For $J \subseteq I$ we write $\Omega_J := (\omega_i)_{i \in J}$. A family $\Omega_J$ with $J \subseteq I$ is a **basis of $\Omega$ along $N$**, if $\Omega_J$ is transverse to $N$ and

$$T_y N \cap \bigcap_{i \in I} \ker(\omega_i(y)) = T_y N \cap \bigcap_{i \in J} \ker(\omega_i(y))$$

for all $y \in N$, or equivalently, if for each $y \in N$ the family $\bigl(\omega_i(y)\big|T_yN\bigr)_{i \in J}$ is a basis of the linear subspace of $(T_yN)^*$ generated by $\bigl(\omega_i(y)\big|T_yN\bigr)_{i \in I}$.

**1.6. Lemma.** *Let $\Omega = (\omega_i)_{i \in I}$ be a finite family of 1-forms of class $C^1$ on an open set $U \subseteq \mathbb{R}^n$, and let $N$ be a submanifold of $U$. Let $J \subseteq I$ and suppose that $\Omega_J$ is a basis of $\Omega$ along $N$. For each $i \in I$ let $L_i$ be a Rolle leaf of $\omega_i = 0$, and write $W_I := N \cap \bigcap_{i \in I} L_i$ and $W_J := N \cap \bigcap_{i \in J} L_i$.*

*Then $W_J$ is either empty or a manifold of dimension $\dim(N) - |J|$, and each component of $W_I$ is a component of $W_J$. In particular, each component of $W_I$ is a manifold of dimension $\dim(N) - |J|$.*

*Proof.* Assume that $W_J$ is not empty. We first show that $W_J$ is a manifold of dimension $\dim(N) - |J|$. We may clearly assume that $|J| = 1$, and write $L$ for the single corresponding Rolle leaf. Fixing $x \in N \cap L$, and working locally around $x$ (see for instance [1]), we further reduce to the case that $x = 0$ and that there is an open set $V \subseteq \mathbb{R}^n$ such that $0 \in V$ and $V \cap L = V \cap (\{0\} \times \mathbb{R}^{n-1})$; it is then enough to show that $V \cap N \cap L$ is a manifold of dimension $\dim(N) - 1$.

The transversality hypothesis now implies that $\Pi_1 | (V \cap N)$ has rank 1 at 0, where $\Pi_1 : \mathbb{R}^n \longrightarrow \mathbb{R}$ is the projection on the first coordinate. Therefore, after shrinking $V$ if necessary, $\Pi_1 | (V \cap N)$ has constant rank 1, and so by the Rank Theorem [13] each nonempty fiber $(V \cap N)_{x_1}$ with $x_1 \in \mathbb{R}$ is a manifold of dimension $\dim(N) - 1$. Therefore $V \cap N \cap L = \{0\} \times (V \cap N)_0$ is a manifold of dimension $\dim(N) - 1$. This finishes the proof of the first assertion.

Let now $C$ be a component of $W_J$. Then $C$ is a manifold of dimension $\dim(N) - |J|$, because $C$ is an open subset of the manifold $W_J$. Let $L := L_j$ for some $j \in I \setminus J$ be such that $C \cap L \neq \emptyset$. Since $\Omega_J$ is a basis of $\Omega$ along $N$, we get for all $y \in C$ that
$$T_yC = T_yN \cap \bigcap_{i \in J} \ker(\omega_i(y)) \subseteq \ker(\omega_j(y)).$$
Therefore, Lemma 1.4 implies that $C \subseteq L$. Since $C$ is arbitrary, this shows that every component of $W_I$ is a component of $W_J$. $\square$

**1.7. Remark.** Let $\omega$ be a 1-form of class $C^1$ on an open set $U \subseteq \mathbb{R}^n$, and let $\sigma = (\sigma_1, \ldots, \sigma_n) : V \longrightarrow U$ be a $C^1$ map, where $V \subseteq \mathbb{R}^m$ is open. Recall [13] that the **pullback** $\sigma^*\omega$ is the 1-form (of class $C^1$) on $V$ given by
$$\sigma^*\omega = (a_1 \circ \sigma)d\sigma_1 + \cdots + (a_n \circ \sigma)d\sigma_n$$
$$= \sum_{i=1}^m \left((a_1 \circ \sigma)\frac{\partial \sigma_1}{\partial y_i} + \cdots + (a_n \circ \sigma)\frac{\partial \sigma_n}{\partial y_i}\right) dy_i.$$

Assume in addition that $\sigma$ is a diffeomorphism onto a submanifold $N$ of $U$. If $\omega$ is transverse to $N$ and $L \subseteq U$ is a Rolle leaf of $\omega = 0$, then $\sigma^*\omega$ is nonsingular and each component of $\sigma^{-1}(L \cap N)$ is a Rolle leaf of $\sigma^*\omega = 0$. More generally, if $\Omega = (\omega_i)_{i \in I}$ is a finite family of 1-forms on $U$ which is transverse to $N$, then the family $\sigma^*\Omega := (\sigma^*\omega_i)_{i \in I}$ is transverse (to $W$).

## 2. Khovanskii Theory over an o-Minimal Structure

Van den Dries noticed some time ago that the Khovanskii theory in [10] could be adapted to the o-minimal setting. In particular, he proved Lemma 2.5 below, replacing a rather special semianalytic argument about carpeting functions in [10].

From now on let $\widetilde{\mathbb{R}}$ be an o-minimal expansion of the real field. (See also [2, 3, 12] for general background on o-minimal structures.) Throughout this section "definable" means "definable in $\widetilde{\mathbb{R}}$ with parameters in $\mathbb{R}$", and $U \subseteq \mathbb{R}^n$ is a definable open set and $\Omega = (\omega_1, \ldots, \omega_q)$ a finite family of **definable 1-forms** of class $C^1$ on $U$, that is, each $\omega_i$ is of the form
$$\omega_i = a_{i,1}dx_1 + \cdots + a_{i,n}dx_n$$
with each $a_{i,j} : U \longrightarrow \mathbb{R}$ a definable $C^1$ function.

**2.1. Lemma.** *Let $A \subseteq \mathbb{R}^n$ be a definable set. Then there is a finite partition $\mathcal{P}$ of $U$ into definable $C^1$ cells such that $\mathcal{P}$ is compatible with $A$, and for each $N \in \mathcal{P}$ there exists a $J \subseteq \{1, \ldots, q\}$ such that $\Omega_J$ is a basis of $\Omega$ along $N$.*

*Proof.* By induction on $d := \dim(A)$; the case $d = 0$ is trivial. So assume that $d > 0$ and that the lemma holds for lower values of $d$. By $C^1$ cell decomposition and the inductive hypothesis, we reduce to the case that $A = N$ is a $C^1$ cell of dimension $d$ (and hence in particular a manifold).

Now for $x \in N$ and $J \subseteq \{1, \ldots, q\}$ we write $T_x(\Omega_J) := T_x N \cap \bigcap_{i \in J} \ker(\omega_i(x))$. For each $J \subseteq \{1, \ldots, q\}$ we define the set
$$N_J := \{ x \in N : T_x(\Omega_J) = T_x(\Omega) \text{ and } \dim(T_x(\Omega_J)) = d - |J| \}.$$

Clearly the sets $N_J$ for $J \subseteq \{1, \ldots, q\}$ form a covering of $N$, and since each form $\omega_i$ is definable, each set $N_J$ is definable. Let $\mathcal{D}$ be a $C^1$ cell decomposition compatible with $U$, $N$, and each $N_J$ for $J \subseteq \{1, \ldots, q\}$. If $C \in \mathcal{D}$ with $C \subseteq N$ and $\dim(C) = d$, and if $J \subseteq \{1, \ldots, q\}$ is such that $C \subseteq N_J$, then $\Omega_J$ is a basis of $\Omega$ along $C$ (since $C$ is open in $N$). On the other hand, let $B := \bigcup \{ C \in \mathcal{D} : \dim(C) < d \}$; then $B$ is definable and $\dim(B) < d$, so the lemma holds with $B$ in place of $A$ by the inductive hypothesis. $\square$

**2.2. Proposition.** *Let $N$ be a definable $C^1$ cell contained in $U$, and suppose that $d := \dim(N) > q$ and that $\Omega$ is transverse to $N$. Then there is a definable closed subset $B$ of $N$ with $\dim(B) < d$, such that whenever $L_i$ is a Rolle leaf of $\omega_i = 0$ for each $i$, then each component of $N \cap L_1 \cap \cdots \cap L_q$ meets $B$.*

Some extra notation is needed to prove Proposition 2.2.

**2.3. Definition.** Given a cell $C \subseteq \mathbb{R}^n$, a **positive form** on $C$ is a definable proper map $\phi : C \longrightarrow \mathbb{R}$ such that $\phi(x) > 0$ for all $x \in C$; or equivalently, it is a definable continuous function $\phi : C \longrightarrow (0, \infty)$ such that
  (i) for each $y \in \mathrm{fr}(C)$ we have $\lim_{x \to y} \phi(x) = +\infty$;
  (ii) for each positive real $r$ the set $\phi^{-1}([0, r])$ is bounded in $\mathbb{R}^n$.
(Of course, (ii) is automatic if $C$ is itself bounded.)

**2.4. Example.** Let $u_1, \ldots, u_n$ be positive real numbers. Then the quadratic form
$$x \mapsto \phi_u(x) := u_1 x_1^2 + \cdots + u_n x_n^2$$
is a positive form on $\mathbb{R}^n$. More generally, if $m \geq n$ and $C \subseteq \mathbb{R}^m$ is a definable cell of dimension $n$, then by [3] there is a definable homeomorphism $\sigma : C \longrightarrow \mathbb{R}^n$. Thus, the function $\phi_u \circ \sigma : C \longrightarrow \mathbb{R}$ is a positive form on $C$.

**2.5. Lemma.** *Let $N$ be a definable $C^1$ cell contained in $U$, and suppose that $d := \dim(N) > q$ and that $\Omega$ is transverse to $N$. Then there is a positive form $\phi$ on $N$ of class $C^1$ such that the (definable) set*
$$Y := \left\{\, a \in N : \, d\phi(a)\big|T_a N \in \langle\, \omega_1(a)\big|T_a N, \, \ldots, \, \omega_q(a)\big|T_a N \,\rangle \,\right\}$$
*has dimension less than $n$.*

*Proof.* By [3] there is a definable diffeomorphism $\sigma : \mathbb{R}^d \longrightarrow N$ of class $C^1$, where $d := \dim(N)$. Replacing $n$ by $d$, $N$ by $\mathbb{R}^d$ and each $\omega_i$ by its pullback $\sigma^* \omega_i$, we may reduce to the case that $N = \mathbb{R}^n$ and each $\omega_i$ is a definable 1-form of class $C^1$ on $\mathbb{R}^n$ with $S(\omega_i) = \emptyset$.

For $u = (u_1, \ldots, u_n) \in \mathbb{R}^n$ with all $u_i > 0$, let
$$B_u := \left\{\, a \in \mathbb{R}^n : \, d\phi_u(a) \in \langle\, \omega_1(a), \, \ldots, \, \omega_q(a) \,\rangle \,\right\},$$
where $\phi_u$ is as in Example 2.4. If $\dim(B_u) < n$ for some $u$ as above, the proof is finished. So assume for a contradiction that $\dim(B_u) = n$ for all $u$ as above. Then $\dim(B) = 2n$, where
$$B := \{\, (u, a) \in \mathbb{R}^n \times \mathbb{R}^n : \, u_1 > 0, \, \ldots, \, u_n > 0, \, a \in B_u \,\},$$
so there are nonempty open $V \subseteq (0, \infty)^n$ and $W \subseteq \mathbb{R}^n$ such that $V \times W \subseteq B$. Fix some $a \in W$ with all $a_i \neq 0$ and let $u$ vary over $V$. Note that
$$d\phi_u(a) \in \langle\, \omega_1(a), \, \ldots, \, \omega_q(a) \,\rangle \iff d\phi_u(a) \wedge \omega_1(a) \wedge \cdots \wedge \omega_q(a) = 0 \quad (\text{in } \bigwedge^{q+1} \mathbb{R}^n).$$
Now $d\phi_u(a) = 2u_1 a_1 dx_1(a) + \cdots + 2u_n a_n dx_n(a)$, so
$$d\phi_u(a) \wedge \omega_1(a) \wedge \cdots \wedge \omega_q(a)$$
$$= 2u_1 a_1 \big(\, dx_1(a) \wedge \omega_1(a) \wedge \cdots \wedge \omega_q(a) \,\big) + \cdots + 2u_n a_n \big(\, dx_n(a) \wedge \omega_1(a) \wedge \cdots \wedge \omega_q(a) \,\big).$$
Since $n > q$, at least one of the covectors $dx_i(a) \wedge \omega_1(a) \wedge \cdots \wedge \omega_q(a)$ is nonzero, and it follows that we can choose a $u \in V$ such that $d\phi_u(a) \wedge \omega_1(a) \wedge \cdots \wedge \omega_q(a) \neq 0$, i.e. $a \notin B_u$, a contradiction. $\square$

*Proof of Proposition 2.2.* Apply the previous lemma to obtain a positive form $\phi$ on $N$ of class $C^1$ and a corresponding definable set $Y \subseteq N$ of dimension less than $d$. This set $Y$ has the desired property: let $C$ be a component of $N \cap L_1 \cap \cdots \cap L_q$, where each $L_i$ is a Rolle leaf of $\omega_i = 0$. Then $C$ is closed in $N$, since each $L_i$ is closed in $U \setminus S(\omega_i)$. Thus $\phi|C$ assumes a minimum value, say at the point $x$ in $C$. This means precisely that $d\phi(x)|T_x C \in \langle\, \omega_1(x)|T_x C, \, \ldots, \, \omega_q(x)|T_x C \,\rangle$, that is, $x \in Y$. $\square$

**2.6. Theorem.** *Let $A \subseteq \mathbb{R}^n$ be a definable set. Then there exists a $K \in \mathbb{N}$ such that whenever $L_i$ is a Rolle leaf of $\omega_i = 0$ for each $i$, then $A \cap L_1 \cap \cdots \cap L_q$ is a union of at most $K$ connected manifolds.*

*Proof.* We proceed by induction on $d := \dim(A)$ and $q$. The cases $d = 0$ or $q = 0$ being trivial, we assume that $d > 0$ and $q > 0$ and that the result holds for lower values of $d$ and $q$. By Lemmas 1.6 and 2.1, it suffices to consider the case that $A = N$ is a $C^1$ cell contained in $U$ and $\Omega$ is transverse to $N$. Note that then $d \geq q$. For each $i$ we let $L_i$ be a Rolle leaf of $\omega_i = 0$, and we put $L := L_1 \cap \cdots \cap L_q$.

**Case $d = q$.** Let $\Omega' := (\omega_1, \ldots, \omega_{q-1})$ and put $L' := L_1 \cap \cdots \cap L_{q-1}$. Then $\dim(N \cap L') = 1$. By the inductive assumption there is a $K \in \mathbb{N}$ (depending only on $N$ and $\Omega'$, but not on the particular Rolle leaves) such that the manifold $N \cap L'$ has at most $K$ components. Let $C$ be a component of $N \cap L'$. If $C \cap L_q$ has more than one point, then by the Rolle property of $L_q$ (and the fact that $C$ is a connected $C^1$ manifold of dimension 1), $C$ is tangent at some point $x \in C$ to the hyperplane field defined by $\omega_q = 0$, which contradicts that $\omega_q$ is transverse to $N \cap L'$. So $C \cap L_q$ has at most one point, for each component $C$ of $N \cap L'$. Hence $|N \cap L| \leq K$.

**Case $d > q$.** Let $Y$ be a closed definable subset of $N$ with the property described in Proposition 2.2, in particular $\dim(Y) < d$. By the inductive hypothesis there is a $K \in \mathbb{N}$, independent of the particular Rolle leaves chosen, such that $Y \cap L$ has at most $K$ components. Hence by Proposition 2.2, the set $N \cap L$ has at most $K$ components. Lemma 1.6 now implies that each component of $L \cap N$ is a manifold, so the theorem is proved. □

**2.7. Corollary.** *Assume that $1 \leq m \leq n$ and let $A \subseteq \mathbb{R}^n$ be a definable set. Then there is a $K \in \mathbb{N}$ such that whenever $a \in \mathbb{R}^m$ and $L_i$ is a Rolle leaf of $\omega_i = 0$ for each $i$, then the fiber $(A \cap L_1 \cap \cdots \cap L_q)_a$ is a union of at most $K$ connected manifolds.*

*Proof.* Note that for each real constant $c$ and each $l \in \{1, \ldots, m\}$, the subset of $U$ defined by the equation $x_l = c$ is a Rolle leaf of the definable 1-form $dx_l = 0$. Applying Theorem 2.6 with $A$ and $\widetilde{\Omega} := (\omega_1, \ldots, \omega_q, dx_1, \ldots, dx_m)$ in place of $\Omega$, we obtain a bound $K \in \mathbb{N}$ such that whenever $a \in \mathbb{R}^m$ and $L_i$ is a Rolle leaf of $\omega_i = 0$ for each $i = 1, \ldots q$, then the fiber

$$(A \cap L_1 \cap \cdots \cap L_q)_a = \Pi\big(A \cap L_1 \cap \cdots \cap L_q \cap \{x \in U : x_1 = a_1\} \cap \cdots \cap \{x \in U : x_m = a_m\}\big)$$

is a union of at most $K$ connected manifolds (here $\Pi : \mathbb{R}^n \longrightarrow \mathbb{R}^{n-m}$ denotes the projection on the last $n - m$ coordinates). □

In the next section the following improvement over Lemma 2.1 will be used.

**2.8. Lemma.** *Let $A \subseteq \mathbb{R}^n$ be a definable set. Then there is a finite partition $\mathcal{P}$ of $U$ into definable $C^1$ cells, such that $\mathcal{P}$ is compatible with $A$ and for each $N \in \mathcal{P}$ and each $J \subseteq \{1, \ldots, q\}$ there exists a $J' \subseteq J$ such that $\Omega_{J'}$ is a basis of $\Omega_J$ along $N$.*

*Proof.* By induction on $d := \dim(A)$; the case $d = 0$ is again trivial. So assume $d > 0$ and that the lemma holds for lower values of $d$. We apply Lemma 2.1 with each $\Omega_J$ in place of $\Omega$ for $J \subseteq \{1, \ldots q\}$ to obtain a corresponding $C^1$ cell decomposition $\mathcal{D}_J$. Let $\mathcal{D}$ be a $C^1$ cell decomposition that is a common refinement of all $\mathcal{D}_J$.

If $C \in \mathcal{D}$ with $\dim(C) = d$, then for every $J \subseteq \{1, \ldots, q\}$ and every $D \in \mathcal{D}_J$ such that $C \subseteq D \subseteq A$, we have that $\dim(D) = \dim(C)$, so that $C$ is open in $D$. Therefore, for every $J \subseteq \{1, \ldots, q\}$ there is a subset $J' \subseteq J$ such that $\Omega_{J'}$ is a basis of $\Omega_J$ along $C$. On the other hand, the set $B := \bigcup \{ C \in \mathcal{D} : \dim(C) < d \}$ satisfies $\dim(B) < d$, so the proof is complete by the inductive hypothesis. $\square$

**2.9. Definition.** A manifold $M \subseteq \mathbb{R}^n$ is **in standard position** if for every strictly increasing map $\iota : \{1, \ldots, k\} \longrightarrow \{1, \ldots, n\}$ there is a $d = d(M, \iota) \le \dim(M)$ such that $\Pi_\iota | M$ has constant rank $d$.

*Remark.* Let $M \subseteq \mathbb{R}^n$ be a manifold in standard position. Then by the rank theorem, for every strictly increasing map $\iota : \{1, \ldots, k\} \longrightarrow \{1, \ldots, n\}$ there is an $e = e(M, \iota) \le \dim(M)$ such that for every $a \in \mathbb{R}^k$ the fiber

$$M_{\iota, a} := \Pi_{\widehat{\iota}}\left(M \cap \Pi_\iota^{-1}(a)\right) = \{ y \in \mathbb{R}^{n-k} : \exists x \in M\big(\Pi_\iota(x) = a \text{ and } \Pi_{\widehat{\iota}}(x) = y\big) \}$$

is either empty or a manifold of dimension $e$, where $\widehat{\iota} : \{1, \ldots, n-k\} \longrightarrow \{1, \ldots, n\}$ is the unique strictly increasing map satisfying $\iota(\{1, \ldots, k\}) \cup \widehat{\iota}(\{1, \ldots, n-k\}) = \{1, \ldots, n\}$.

**2.10. Definition.** A manifold $N \subseteq U$ is **in $\Omega$-position** if there is a $J \subseteq \{1, \ldots, q\}$ such that $\Omega_J$ is a basis of $\Omega$ along $N$, and whenever $L_i$ is a Rolle leaf of $\omega_i = 0$ for each $i$, then $N \cap L_1 \cap \cdots \cap L_q$ is in standard position.

**2.11. Corollary.** *Let $A \subseteq \mathbb{R}^n$ be a definable set. Then there is a finite partition $\mathcal{P}$ of $U$ into definable $C^1$ cells, such that $\mathcal{P}$ is compatible with $A$ and each $N \in \mathcal{P}$ is in $\Omega$-position.*

*Proof.* Apply Lemma 2.8 with $\widetilde{\Omega} := (\omega_1, \ldots, \omega_q, dx_1, \ldots, dx_n)$ in place of $\Omega$ to obtain a corresponding partition $\mathcal{P}$. Then, given $N \in \mathcal{P}$ there is a $J \subseteq \{1, \ldots, q\}$ such that $\Omega_J$ is a basis of $\Omega$ along $N$ (since $\Omega \subseteq \widetilde{\Omega}$). In particular, given any Rolle leaves $L_i$ of $\omega_i = 0$ for each $i$, the set $N \cap L_1 \cap \cdots \cap L_q$ is either empty or a manifold of dimension $\dim(N) - |J|$.

It remains to show that $N \cap L_1 \cap \cdots \cap L_q$ is in standard position. Let $\iota : \{1, \ldots, k\} \longrightarrow \{1, \ldots, n\}$ be a strictly increasing map, and let $a \in \mathbb{R}^k$. Since the hyperplane defined by the equation $x_{\iota(j)} = a_j$ is a Rolle leaf of $dx_{\iota(j)} = 0$ for each $j = 1, \ldots, k$, the above application of Lemma 2.8 implies that the set $N(a) := N \cap L_1 \cap \cdots \cap L_q \cap \{x_{\iota(1)} = a_1\} \cap \cdots \cap \{x_{\iota(k)} = a_k\}$ is either empty or a manifold of dimension $d$, with $d$ only depending on $N$ and $\iota$. Moreover, for any $x \in N(a)$ we have

$$\ker\left(\Pi_\iota \big| T_x(N \cap L_1 \cap \cdots \cap L_q)\right) = T_x(N \cap L_1 \cap \cdots \cap L_q) \cap \ker(\Pi_\iota)$$
$$= T_x N \cap \bigcap_{i=1}^{q} \ker(\omega_i(x)) \cap \bigcap_{j=1}^{k} \ker(dx_{\iota(j)}(x))$$
$$= T_x N(a).$$

Hence $\ker\left(\Pi_\iota \big| T_x(N \cap L_1 \cap \cdots \cap L_q)\right)$ has dimension $d$, and since $d$ is independent of $x$ and $a$, this shows that $\Pi_\iota \big|(N \cap L_1 \cap \cdots \cap L_q)$ has constant rank. $\square$

## 3. $T^\infty$-Sets

We define in this section the notion of a "$T^\infty$-set", which Lion and Rolin [8] introduced in the analytic context (under the name "$T^\infty$-Pfaffian set"). We show that every $T^\infty$-set has finitely many connected components and that the collection of $T^\infty$-subsets of $I^n$ (for various $n$) forms a structure on $I := [0, 1]$. (In contrast to [8], we do not use the "Cauchy-Crofton formula" here.) In particular, this structure is o-minimal.

We fix an o-minimal expansion $\widetilde{\mathbb{R}}$ of the real field, and we let $i$ and $j$ range over $\mathbb{N}$. As before, "definable" means "definable in $\widetilde{\mathbb{R}}$".

**3.1. Definition.** A set $W \subseteq \mathbb{R}^n$ is called **Pfaffian** if there are definable 1-forms $\omega_1, \ldots, \omega_q$ of class $C^1$ on a definable open set $U \subseteq \mathbb{R}^n$, Rolle leaves $L_p$ of $\omega_p = 0$ for each $p = 1, \ldots, q$, and a definable set $A \subseteq U$ such that
$$W = A \cap L_1 \cap \cdots \cap L_q.$$

*Remark.* Let $W \subseteq \mathbb{R}^n$ be Pfaffian.
  (1) If $V \subseteq \mathbb{R}^n$ is Pfaffian, then $V \cap W$ is Pfaffian.
  (2) If $n' \in \mathbb{N}$ and $W' \subseteq \mathbb{R}^{n'}$ is Pfaffian, then $W \times W' \subseteq \mathbb{R}^{n+n'}$ is Pfaffian.
  (3) If $\pi$ is a permutation of $\{1, \ldots, n\}$, then the set $\pi(W) := \{\, (x_{\pi(1)}, \ldots, x_{\pi(n)}) : x \in W \,\}$ is Pfaffian.
  (4) If $n' \leq n$, then for every $a \in \mathbb{R}^{n'}$ the fiber $W_a \subseteq \mathbb{R}^{n-n'}$ is Pfaffian.

**3.2. Definition.** A set $X \subseteq \mathbb{R}^m$ is a **basic $T^\infty$-set** if there are $k$, $l$, parameters $\epsilon(i, j) \in \mathbb{R}^k$ for each $i$ and $j$, and a Pfaffian set $W \subseteq \mathbb{R}^k \times \mathbb{R}^m \times \mathbb{R}^l$, such that
  (i) for each pair $(i, j)$ the fiber $W_{\epsilon(i,j)} \subseteq \mathbb{R}^m \times \mathbb{R}^l$ is compact;
  (ii) for each $i$ the sequence $\big(W(i, j)\big)_j$ of subsets of $\mathbb{R}^m$ is decreasing, where $W(i, j) := \Pi(W_{\epsilon(i,j)})$ and $\Pi : \mathbb{R}^{m+l} \longrightarrow \mathbb{R}^m$ is the projection on the first $m$ coordinates;
  (iii) the sequence $\big(W(i)\big)_i$ is increasing, where $W(i) := \bigcap_j W(i, j)$; and
  (iv) $X = \bigcup_i W(i)$.

In this case we say that $X$ is **obtained from** $W$. A finite union of basic $T^\infty$-sets is a **$T^\infty$-set**.

*Remark.* Let $X \subseteq \mathbb{R}^m$ be a basic $T^\infty$-set obtained from $W \subseteq \mathbb{R}^n$.
  (1) If $\pi$ is a permutation of $\{1, \ldots, m\}$, then $\pi(X) := \{\, (x_{\pi(1)}, \ldots, x_{\pi(m)}) : x \in X \,\}$ is a basic $T^\infty$-set.
  (2) If $m' \leq m$, then for every $a \in \mathbb{R}^{m'}$ the fiber $X_a$ is a basic $T^\infty$-set obtained from $W$ (namely, if $\epsilon(i, j) \in \mathbb{R}^k$ are the parameters involved in the definition of $X$, then the parameters for $X_a$ are $(\epsilon(i, j), a) \in \mathbb{R}^{k+m'}$).
  (3) If $m' \leq m$, then the projection $\Pi_{m'}(X)$ of $X$ on the first $m'$ coordinates is a basic $T^\infty$-set obtained from $W$ (replacing $m$ and $l$ in the definition of $X$ by $m'$ and $m - m' + l$ for $\Pi_{m'}(X)$).

**3.3. Lemma.** *Every Pfaffian set is a $T^\infty$-set.*

*Proof.* By $C^1$ cell decomposition it suffices to consider $W = A \cap L_1 \cap \cdots \cap L_q \subseteq \mathbb{R}^n$ as in Definition 3.1, where $A \subseteq U$ is a definable $C^1$ cell. Let $g$ be a positive form on $A$ (as given by Example 2.4, say), and put $\widetilde{W} := \widetilde{A} \cap (\mathbb{R} \times L_1) \cap \cdots \cap (\mathbb{R} \times L_q)$ with

$$\widetilde{A} := \{ (\theta, x) \in (0, \infty) \times A : g(x) \leq 1/\theta \} \subseteq \mathbb{R}^{1+n}.$$

By Remark 1.7 each $\mathbb{R} \times L_p$ is a Rolle leaf of the pullback $\Pi^* \omega_p$ on $\mathbb{R} \times U$, where $\Pi : \mathbb{R}^{1+n} \longrightarrow \mathbb{R}^n$ is the projection on the last $n$ coordinates. Since each $\mathbb{R} \times L_i$ is a closed subset of $\mathbb{R} \times U$, each fiber $\widetilde{W}_\theta$ with $\theta > 0$ is compact. Now $W$ can be written as the increasing union of the compact sets $\widetilde{W}_{1/i}$ for nonzero $i \in \mathbb{N}$. □

**3.4. Proposition.** *Let $W \subseteq \mathbb{R}^n$ be a Pfaffian set. Then there is an $N \in \mathbb{N}$ such that every basic $T^\infty$-set obtained from $W$ has at most $N$ components. In particular, every $T^\infty$-set has finitely many components.*

*Proof.* From Corollary 2.7 we obtain an $N \in \mathbb{N}$ which bounds the number of components of any fiber of $W$. Next note that if $K \subseteq \mathbb{R}^m$ is the intersection (resp. union) of a decreasing (resp. increasing) sequence $\big(K(i)\big)_i$ of compact subsets of $\mathbb{R}^m$ and if each $K(i)$ has at most $N$ components, then $K$ also has at most $N$ components. □

**3.5. Proposition.** *The collection of $T^\infty$-sets is closed under taking finite unions, finite intersections, projections and topological closure.*

*Proof.* Closure under taking finite unions and projections is obvious from the definition of $T^\infty$-set. For finite intersections, let $X, X' \subseteq \mathbb{R}^m$ be basic $T^\infty$-sets obtained from $W \subseteq \mathbb{R}^n$ and $W' \subseteq \mathbb{R}^{n'}$, respectively. Then $X \times X'$ is a basic $T^\infty$-set obtained from $\widetilde{W} := \{ (\theta, \theta', x, x', y, y') \in \mathbb{R}^{n+n'} : (\theta, x, y) \in W, (\theta', x', y') \in W' \}$. Therefore the set $(X \times X') \cap \Delta$ is a basic $T^\infty$-set obtained from $\widetilde{W} \cap \widetilde{\Delta}$, where $\Delta := \{ (x, x') \in \mathbb{R}^m \times \mathbb{R}^m : x = x' \}$ and $\widetilde{\Delta} := \{ (\theta, \theta', x, x', y, y') \in \mathbb{R}^{n+n'} : x = x' \}$. But $X \cap X' = \Pi_m\big((X \times X') \cap \Delta\big)$, which finishes the proof for finite intersections.

For topological closure, let $X \subseteq \mathbb{R}^m$ be a basic $T^\infty$-set obtained from $W \subseteq \mathbb{R}^n$ as in 3.2, with $n = k + m + l$. We define the Pfaffian set $\widetilde{W} \subseteq \mathbb{R}^{\widetilde{n}}$ with $\widetilde{n} := k + 1 + m + l + m + l$ by

$$\widetilde{W} := \{ (\theta, \eta, x, y, x', y') : (\theta, x', y') \in W, |(x, y) - (x', y')| \leq \eta \}.$$

Note that $\widetilde{\Pi}(\widetilde{W}_{\epsilon(i,j), \eta}) = \mathrm{cl}\big(T(W(i, j), \eta)\big)$ for all $i, j$ and $\eta$, where $\widetilde{\Pi} : \mathbb{R}^{m+l+m+l} \longrightarrow \mathbb{R}^m$ is the projection on the first $m$ coordinates. We claim that there are sequences $\big(\theta(i)\big)_i$ with each $\theta(i) \in \mathbb{R}^k$ and $\big(\eta(i)\big)_i$ with each $\eta(i) \in \mathbb{R}$, such that $\big(\widetilde{\Pi}\big(\widetilde{W}_{\theta(i), \eta(i)}\big)\big)_i$ is a decreasing sequence of compact sets with intersection $\mathrm{cl}(X)$, which then finishes the proof.

We use the notation introduced in Definition 3.2. Clearly each $\widetilde{W}_{\epsilon(i,j), \eta}$ is compact. Let $\delta(i, j) := \sup_{x \in W(i,j)} d(x, W(i))$ and $\delta(i) := \sup_{x \in \mathrm{cl}(X)} d(x, W(i))$. For each $i$ the sequence $\big(\delta(i, j)\big)_j$ converges to 0, and the sequence $\big(\delta(i)\big)_i$ converges to 0. If $\delta(i) = 0$ for some $i$ the claim is trivial, so we assume that each $\delta(i)$ is nonzero. Passing to a subsequence if necessary, we may also assume that $\delta(i+1) \leq \delta(i)/3$. For each $i$ we fix a $j_i$ such that $\delta(i, j_i) \leq \delta(i)$. Then

$$T\big(\mathrm{cl}(X), \delta(i)\big) \subseteq T\big(W(i, j_i), 2\delta(i)\big) \subseteq T\big(\mathrm{cl}(X), 3\delta(i)\big),$$

so $\delta(i+1) \leq \delta(i)/3$ implies that $T\big(W(i+1, j_{i+1}), 2\delta(i+1)\big) \subseteq T\big(W(i, j_i), 2\delta(i)\big)$. Therefore, $\theta(i) := \epsilon(i, j_i)$ and $\eta(i) := 2\delta(i)$ will do. □

**3.6. Remark.** Let $X, X' \subseteq \mathbb{R}^m$ be basic T$^\infty$-sets obtained from the Pfaffian sets $W$ and $W'$ respectively. The proof above shows that $X \cap X'$ is a basic T$^\infty$-set obtained from a Pfaffian set $\widetilde{W}$ that depends only on $W$ and $W'$ (but not on the particular sets $X$ and $X'$).

The difficult part is to prove that the collection of T$^\infty$-subsets of $[0,1]^m$ is closed under taking complements. The main step in the proof is to show that the boundary of a bounded T$^\infty$-set is contained in a closed T$^\infty$-set with empty interior (see Lemma 3.11 below). We will do this by induction on $m$, using a fibering argument. First we show that for a T$^\infty$-set, taking the closure of its fibers is "almost everywhere" the same as taking the fibers of its closure.

**3.7. Lemma.** *Let $X \subseteq \mathbb{R}^m$ be a T$^\infty$-set with $m \geq 1$. Then the set*
$$B := \{\, a \in \mathbb{R} : \operatorname{cl}(X_a) \neq \operatorname{cl}(X)_a \,\}$$
*is countable.*

*Proof.* The following short proof is due to Chris Miller. The case $m = 1$ follows from Proposition 3.4, so we assume that $m > 1$. For each $a \in B$ there is a box $U \subseteq \mathbb{R}^{m-1}$ such that $\operatorname{cl}(X_a) \cap U = \emptyset$, but $\operatorname{cl}(X)_a \cap U \neq \emptyset$. Hence $B = \bigcup_U B_U$, where $U$ ranges over all rational boxes in $\mathbb{R}^{m-1}$ and
$$B_U := \{\, a \in \mathbb{R} : \operatorname{cl}(X_a) \cap U = \emptyset,\ \operatorname{cl}(X)_a \cap U \neq \emptyset \,\}.$$
One easily verifies that for each $U$ the set $B_U$ is contained in the frontier of the T$^\infty$-set $\Pi(X \cap (\mathbb{R} \times U))$, where $\Pi : \mathbb{R}^m \longrightarrow \mathbb{R}$ is the projection on the first coordinate. So by Proposition 3.4 each $B_U$ is finite. □

Next we prove a lemma similar to the above, but with "Hausdorff limit" in place of "closure" (see Lemma 3.9 below). Recall that for any two nonempty compact sets $S, T \subseteq \mathbb{R}^m$ the **Hausdorff distance** $d(S, T)$ is the greater of the two values $\max\{d(x, T) : x \in S\}$ and $\max\{d(y, S) : y \in T\}$. The set $\mathbb{H}_m$ of all nonempty compact subsets of $\mathbb{R}^m$ equipped with the Hausdorff metric is a metric space in which every closed and bounded subset of $\mathbb{H}_m$ is compact.

Let $K(i) \subseteq \mathbb{R}^m$ for $i \in \mathbb{N}$ be compact sets, and let $K \subseteq \mathbb{R}^m$ be a nonempty compact set. The sequence $\big(K(i)\big)_i$ **converges** to $K$ if $K(i) \neq \emptyset$ for all $i$, and for every $\epsilon > 0$ there is an $i$ such that $d(K(j), K) < \epsilon$ whenever $j \geq i$. The set $K$ is a **Hausdorff limit** of $\big(K(i)\big)$ if a subsequence of $\big(K(i)\big)$ converges to $K$, or equivalently, if for every $\epsilon > 0$ there are infinitely many $i$ such that $K(i) \neq \emptyset$ and $d(K(i), K) < \epsilon$.

Assume now that the sequence $\big(K(i)\big)$ is **bounded**, that is, there exists an $R > 0$ such that $K(i) \subseteq [-R, R]^m$ for all $i$. If $K(i) \neq \emptyset$ for infinitely many $i$, then $\big(K(i)\big)$ has a (nonempty) Hausdorff limit. If each $K(i)$ is nonempty and $\big(K(i)\big)$ converges to $K$, then $K$ is precisely the set of all accumulation points of sequences $(x_i)$ with each $x_i \in K(i)$ (recall that an **accumulation point** of $(x_i)$ is by definition the limit of a convergent subsequence of $(x_i)$.).

**3.8. Lemma.** *Let $m \geq 1$, and let $\bigl(K(i)\bigr)_i$ be a bounded sequence of nonempty compact subsets of $\mathbb{R}^m$ converging to a nonempty compact set $K$. Assume that $a \in \mathbb{R}$ is such that $K_a \neq \emptyset$ and $K_a$ is not a Hausdorff limit of $\bigl(K(i)_a\bigr)_i$. Then there are open boxes $U_1, \ldots, U_l \subseteq \mathbb{R}^{m-1}$ such that*

(∗) *$K_a \cap U_j \neq \emptyset$ for each $j = 1, \ldots, l$, and for all sufficiently large $i$ there is a $j$ such that $K(i)_a \cap U_j = \emptyset$.*

*Proof.* If $m = 1$ the lemma is trivial, so we assume that $m > 1$. By the definition of Hausdorff limit, there is an $\epsilon > 0$ such that for all sufficiently large $i$ we have either $K(i)_a = \emptyset$, or $K(i)_a \not\subseteq T(K_a, \epsilon)$, or $K_a \not\subseteq T\bigl(K(i)_a, \epsilon\bigr)$. If $K(i)_a = \emptyset$ for all but finitely many $i$, then (∗) is trivial, so we assume that $K(i)_a \neq \emptyset$ for infinitely many $i$. Passing now to the subsequence of all $K(i)$ with $K(i)_a \neq \emptyset$, we may clearly assume that $K(i)_a \neq \emptyset$ for all $i$.

**Claim.** *$K_a \not\subseteq T\bigl(K(i)_a, \epsilon\bigr)$ for all sufficiently large $i$.*

*Proof.* If $K(i)_a \not\subseteq T(K_a, \epsilon)$ for infinitely many $i$, then there is a strictly increasing sequence $(i_j)_j$ and there are points $x_j \in K(i_j)_a$ for each $j$ such that $d(x_j, K_a) \geq \epsilon$. Since $\bigl(K(i)\bigr)$ is bounded, the sequence $(x_j)$ has an accumulation point $x \in \mathbb{R}^{m-1}$. By the remark before the lemma, since $K$ is the limit of the sequence $\bigl(K(i)\bigr)$, we have $x \in K_a$, which contradicts that $d(x_j, K_a) \geq \epsilon$ for all $j$, so the claim is proved.

By the claim, there is a sequence $(x_i)$ in $K_a$ such that $d(x_i, K(i)_a) \geq \epsilon$ for all sufficiently large $i$. Let $A$ be the set of accumulation points of the sequence $(x_i)$; by the remark before the lemma, we have $A \subseteq K_a$, so $A$ is compact. Choose $y_1, \ldots, y_l \in A$ such that the boxes $U_j := B(y_j, \epsilon/2)$ cover $A$. Clearly $x_i \in U_1 \cup \cdots \cup U_l$ for all sufficiently large $i$. Hence for all sufficiently large $i$ there is a $j \in \{1, \ldots, l\}$ such that $d(y_j, K(i)_a) > \epsilon/2$, that is, $K(i)_a \cap U_j = \emptyset$, which finishes the proof. □

**3.9. Lemma.** *Let $m \geq 1$, let $W \subseteq \mathbb{R}^n$ be a Pfaffian set, and for each $i$ let $K(i) \subseteq \mathbb{R}^m$ be a nonempty compact basic $\mathrm{T}^\infty$-set obtained from $W$. Assume that the sequence $\bigl(K(i)\bigr)_i$ is bounded and converges to a nonempty compact set $K \subseteq \mathbb{R}^m$. Then the set*

$$B := \bigl\{\, a \in \mathbb{R} : \ K_a \neq \emptyset \text{ and } K_a \text{ is not a Hausdorff limit of } \bigl(K(i)_a\bigr) \,\bigr\}$$

*is countable.*

*Proof.* By Proposition 3.4 there is an $N \in \mathbb{N}$ such that each $K(i)$ has at most $N$ components. More generally, if $U \subseteq \mathbb{R}^m$ is any box, then by Remark 3.6 each $K(i) \cap U$ is a basic $\mathrm{T}^\infty$-set obtained from a Pfaffian set $W(U)$ depending only on $U$, but not on $i$. Hence by Proposition 3.4 we can find an $N(U) \in \mathbb{N}$ such that each $K(i) \cap U$ has at most $N(U)$ components.

Let now $a \in B$. By Lemma 3.8 we have $B = \bigcup_U B_U$, where $U$ ranges over all finite tuples $U = (U_1, \ldots, U_l)$ of rational boxes in $\mathbb{R}^{m-1}$ and

$$B_U := \bigl\{\, a \in \mathbb{R} : \ \text{condition } (\ast) \text{ holds for } a \,\bigr\}.$$

Fix $U = (U_1, \ldots, U_l)$ and let $N_1, \ldots, N_l \in \mathbb{N}$ be such that each $K(i) \cap (\mathbb{R} \times U_j)$ has at most $N_j$ components. We show that $|B_U| \leq N := 2N_1 + \cdots + 2N_l$, which clearly finishes the proof of the lemma.

Assume for a contradiction that there are distinct elements $a_1, \ldots, a_{N+1} \in B_U$, and choose $\rho > 0$ such that the intervals $I_p := (a_p - \rho, a_p + \rho)$ for $p = 1, \ldots, N+1$ are disjoint. Choose an $i$ so large that for each $p$ there is a $j = j(p) \in \{1, \ldots, l\}$ with $K(i)_{a_p} \cap U_j = \emptyset$. Then for some $j$ there are distinct $p(1), \ldots, p(2N_j + 1) \in \{1, \ldots, N+1\}$ such that $K(i)_{a_{p(q)}} \cap U_j = \emptyset$ for each $q \in \{1, \ldots, 2N_j + 1\}$. On the other hand, since $K_{a_{p(q)}} \cap U_j \neq \emptyset$ and the sequence $(K(i))$ converges to $K$, we may assume that $i$ is so large that $K(i) \cap (I_{p(q)} \times U_j) \neq \emptyset$ for each $q$. Then $K(i) \cap (\mathbb{R} \times U_j)$ has at least $N_j + 1$ components, contradicting our choice of $N_j$. □

**3.10. Remark.** By definition any $T^\infty$-set is a countable union of compact sets. Therefore:
  (i) If $X(i) \subseteq \mathbb{R}^m$ is a $T^\infty$-set with empty interior, $i \in \mathbb{N}$, then $\bigcup_i X(i)$ has empty interior.
  (ii) If $m \geq 1$, then a $T^\infty$-set $X \subseteq \mathbb{R}^m$ has empty interior if and only if $\{ a \in \mathbb{R} : \operatorname{int}(X_a) \neq \emptyset \}$ has empty interior.

**3.11. Lemma.** *Let $X \subseteq \mathbb{R}^m$ be a bounded $T^\infty$-set. Then $\operatorname{bd}(X)$ is contained in a closed $T^\infty$-set with empty interior.*

*Proof.* We may assume that $X$ is a nonempty basic $T^\infty$-set, say obtained from $W \subseteq \mathbb{R}^n$ as in Definition 3.2; we adopt here the notation set up there. Since $X$ is bounded, we may also assume that $W$ is bounded.

Applying Corollary 2.11 with $\Omega = (\omega_1, \ldots, \omega_q)$, we obtain a partition of $A$ into definable $C^1$ cells $A_1, \ldots, A_M$ in $\Omega$-position. For each $p \in \{1, \ldots, M\}$ we put
$$W^p := A_p \cap L_1 \cap \cdots \cap L_q.$$
Then each nonempty fiber $W^p_\epsilon \subseteq \mathbb{R}^{m+l}$ with $\epsilon \in \mathbb{R}^k$ is a manifold, and there is a $d_p \in \mathbb{N}$, independent of $\epsilon$, such that whenever $W^p_\epsilon$ is nonempty the map $\Pi | W^p_\epsilon$ has constant rank $d_p$. Let $S$ be the set of indices $p \in \{1, \ldots, M\}$ for which $d_p < m$.

We call a set $Y \subseteq \mathbb{R}^m$ an **approximation** of $\operatorname{bd}(X)$ whenever $Y$ is obtained from the sets $W^p$ in the following way:
  (i) for each $p$, $i$ and $j$ we let $K^p(i,j) := \Pi\left(\operatorname{cl}\left(W^p_{\epsilon(i,j)}\right)\right) \subseteq \mathbb{R}^m$.
  (ii) If $p$ and $i$ are such that $K^p(i,j) \neq \emptyset$ for infinitely many $j$, then we let $K^p(i)$ be a Hausdorff limit of $(K^p(i,j))_j$; otherwise we put $K^p(i) := \emptyset$.
  (iii) If $p$ is such that $K^p(i) \neq \emptyset$ for infinitely many $i$, then we let $K^p$ be a Hausdorff limit of $(K^p(i))_i$; otherwise we put $K^p := \emptyset$.
  (iv) We let $Y := \bigcup_{p \in S} K^p$.

(These Hausdorff limits exist because $W$ is bounded.) For the rest of the proof we fix an arbitrary approximation $Y$ of $\operatorname{bd}(X)$, and we use the notation established above. Passing to subsequences if necessary, we will assume that whenever $K^p(i) \neq \emptyset$ the sequence $(K^p(i,j))$ converges to $K^p(i)$, and whenever $K^p \neq \emptyset$ the sequence $(K^p(i))$ converges to $K^p$. The following claims establish the conclusion of the lemma for $Y$.

**Claim 1:** *There is an $\widetilde{n} \geq n$, and for each $p$ there is a bounded Pfaffian set $\widetilde{W^p} \subseteq \mathbb{R}^{\widetilde{n}}$, such that the sets $K^p(i)$ and $K^p$ are basic $T^\infty$-sets obtained from $\widetilde{W^p}$. In particular, $Y$ is a $T^\infty$-set.*

*Proof.* We prove the claim for $K^p$; the case of each $K^p(i)$ is handled similarly. We clearly may assume that $K^p \neq \emptyset$. Let $\delta(i,j) := d(K^p(i,j), K^p(i))$ and $\delta(i) := d(K^p(i), K^p)$. Passing to a subsequence if necessary, we may assume that $\delta(i+1) < \delta(i)/6$ for each $i$. Let $g$ be a positive form on $A_p$ (as given by Example 2.4, say). Let $z = (\theta, t, s, x, y, x', y')$ range over $\mathbb{R}^{\widetilde{n}}$ with $\widetilde{n} := k+1+1+m+l+m+l$, $0 < t < 5 \cdot \delta(0)$ and $0 < s < 1$, and consider the bounded Pfaffian set

$$\widetilde{W}^p := \{\, z : \; (\theta, x', y') \in W^p, \; g(\theta, x', y') \leq 1/s, \; |(x,y) - (x',y')| \leq t \,\}.$$

(The reason for introducing $g$ is to make the fibers of $\widetilde{W}^p$ compact.) We write $\widetilde{\Pi}$ for the projection from $\mathbb{R}^{m+l+m+l}$ on the first $m$ coordinates. For each $i$ choose a $j_i \in \mathbb{N}$ and a $\beta(i) \in (0,1)$ such that $\delta(i, j_i) \leq \delta(i)$ and $K^p(i) \subseteq \widetilde{\Pi}\bigl( \widetilde{W}^p_{\epsilon(i,j_i),\, 2\delta(i),\, \beta(i)} \bigr)$. Then for each $i$,

$$T\bigl(K^p, \delta(i)\bigr) \subseteq \widetilde{\Pi}\left( \widetilde{W}^p_{\epsilon(i,j_i),\, 4\delta(i),\, \beta(i)} \right) \subseteq T\bigl(K^p, 6\delta(i)\bigr).$$

Since $\delta(i+1) \leq \delta(i)/6$, the sequence of sets $\widetilde{\Pi}\bigl( \widetilde{W}^p_{\epsilon(i,j_i),\, 4\delta(i),\, \beta(i)} \bigr)$ is a decreasing sequence of compact sets converging to $K^p$. This proves Claim 1.

**Claim 2:** $Y$ contains $\mathrm{bd}(X)$.

*Proof.* Since $X$ is the union of the increasing sequence of compact sets $(W(i))$, every point $a \in \mathrm{bd}(X)$ is the limit of some sequence $(a_i)_i$ with each $a_i \in \mathrm{bd}(W(i))$. On the other hand, if $Y$ is nonempty, then the sequence of compact sets $\bigl(\bigcup_{p \in S} K^p(i)\bigr)$ converges to $Y$. Hence in order to show that $\mathrm{bd}(X) \subseteq Y$ it suffices to show that $\mathrm{bd}(W(i)) \subseteq \bigcup_{p \in S} K^p(i)$ for each $i$.

Fix an $i$ and let $a \in \mathrm{bd}(W(i))$. Note that if $K^p(i)$ is nonempty, then the sequence $\bigl(K^p(i,j)\bigr)_j$ converges to $K^p(i)$. Therefore, it suffices to show that for every $r > 0$ the box $B(a,r)$ contains a point of $\Pi\bigl(W^p_{\epsilon(i,j)}\bigr)$ for some $p \in S$ and some $j = j(i,r)$.

Fix an $r > 0$ and choose a point $b \in \mathbb{R}^m$ and a $\rho > 0$ such that $\mathrm{cl}(B(b,\rho)) \subseteq B(a,r) \setminus W(i)$. Since $W(i)$ is the intersection of the decreasing sequence of compact sets $\bigl(W(i,j)\bigr)_j$, there is a $j = j(i,r)$ such that $W(i,j) \cap B(b,\rho) = \emptyset$. Choose a point $z \in W(i,j) \cap B(a,r)$. By Proposition 3.4, the intersection of $W(i,j)$ with the segment $[z,b]$ is a finite union of closed segments contained in $B(a,r)$, so $W(i,j) \cap [z,b]$ contains an extremal point $c$ closest to $b$. This point $c$ has to belong to some $\Pi\bigl(W^p_{\epsilon(i,j)}\bigr)$ with $p \in S$, since the sets $\Pi\bigl(W^p_{\epsilon(i,j)}\bigr)$ with $p \notin S$ are open subsets of $\mathbb{R}^m$. But $c \in B(a,r)$, so the claim is proved.

**Claim 3:** $Y$ has empty interior.

*Proof.* We proceed by induction on $m$. If $m = 1$, then $d_p \leq 0$ for each $p \in S$, and since the number of components in $\Pi\bigl(W^p_{\epsilon(i,j)}\bigr)$ is uniformly bounded, this implies that $K^p$ is finite for each $p \in S$. So let $m > 1$ and assume that the claim holds for lower values of $m$.

Let $a \in \mathbb{R}$, and note that $X_a = \bigcup_i \bigcap_j \Pi'\bigl(W_{\epsilon(i,j),\,a}\bigr)$, where $\Pi' : \mathbb{R}^{m-1+l} \longrightarrow \mathbb{R}^{m-1}$ is the projection on the first $m-1$ coordinates. By Corollary 2.11 each fiber $W^p_{\epsilon,a}$ is a manifold, and there is an $e_p \leq m-1$, independent of $(\epsilon, a)$, such that whenever $W^p_{\epsilon,a}$ is nonempty, the projection $\Pi'$ restricted to $W^p_{\epsilon,a}$ has constant rank $e_p$. Put $S' := \{\, p \in \{1, \ldots, M\} : \; e_p < m-1 \,\}$.

In analogy to the above, we call a set $Y' \subseteq \mathbb{R}^{m-1}$ an **approximation** of $\mathrm{bd}(X_a)$ whenever $Y'$ is obtained from the sets $W^p$ in the following way:

(i) for each $p$, $i$ and $j$ we let $(K')^p(i,j) := \Pi\left(\text{cl}\left(W^p_{\epsilon(i,j),a}\right)\right) \subseteq \mathbb{R}^{m-1}$.

(ii) If $p$ and $i$ are such that $(K')^p(i,j) \neq \emptyset$ for infinitely many $j$, then we let $(K')^p(i)$ be a Hausdorff limit of $\left((K')^p(i,j)\right)_j$ ; otherwise we put $(K')^p(i) := \emptyset$.

(iii) If $p$ is such that $(K')^p(i) \neq \emptyset$ for infinitely many $i$, then we let $(K')^p$ be a Hausdorff limit of $\left((K')^p(i)\right)_i$ ; otherwise we put $(K')^p := \emptyset$.

(iv) We let $Y' := \bigcup_{p \in S'}(K')^p$.

Consider now the set
$$H := \{\, a \in \mathbb{R} :\ Y_a \neq \emptyset \text{ and } Y_a \text{ is not an approximation of } \text{bd}(X_a) \,\}.$$

By the inductive hypothesis, for each $a \in \mathbb{R} \setminus H$ the set $Y_a$ has empty interior. So if $H$ has empty interior, then Claim 1 and Remark 3.10(ii) imply that the set $Y$ has empty interior, which then proves Claim 3. We therefore need to prove that $\text{int}(H) = \emptyset$; in fact, we show that $H$ is countable.

By definition the set $H$ is contained in the union of the following sets:
$$G^p(i,j) := \left\{\, a \in \mathbb{R} :\ \Pi'\left(\text{cl}\left(W^p_{\epsilon(i,j),\,a}\right)\right) \neq K^p(i,j)_a \,\right\},$$
$$G^p(i) := \left\{\, a \in \mathbb{R} :\ K^p(i)_a \neq \emptyset \text{ and } K^p(i)_a \text{ is not a Hausdorff limit of } \left(K^p(i,j)_a\right) \,\right\},$$
$$G^p := \left\{\, a \in \mathbb{R} :\ K^p_a \neq \emptyset \text{ and } K^p_a \text{ is not a Hausdorff limit of } \left(K^p(i)_a\right) \,\right\},$$
$$G := \left\{\, a \in \mathbb{R} :\ \cup_{p \in S} K^p_a \neq \cup_{p \in T} K^p_a \,\right\}.$$

It is therefore enough to prove that each of these sets is countable. Since $K^p(i,j)_a = \Pi'\left(\text{cl}\left(W^p_{\epsilon(i,j)}\right)_a\right)$, we have
$$G^p(i,j) \subseteq \left\{\, a \in \mathbb{R} :\ \text{cl}\left(W^p_{\epsilon(i,j),\,a}\right) \neq \text{cl}\left(W^p_{\epsilon(i,j)}\right)_a \,\right\},$$
so each $G^p(i,j)$ is countable by Lemma 3.7. Claim 1 and Lemma 3.9 imply that each $G^p(i)$ and each $G^p$ is countable. It remains to show that $G$ is countable.

Note that $S' \subseteq S$, and let $p \in S \setminus S'$; then we must have $d_p = e_p = m-1$. So from the rank theorem we get $\dim\left(W^p_\epsilon\right) = \dim\left(W^p_{\epsilon,a}\right)$ whenever $W^p_{\epsilon,a} \neq \emptyset$. Since both $W^p_\epsilon$ and $\Pi_1\left(W^p_\epsilon\right)$ have finitely many components (where $\Pi_1 : \mathbb{R}^{m+k} \longrightarrow \mathbb{R}$ is the projection onto the first coordinate), this means that $W^p_{\epsilon,a} = \emptyset$ for all but finitely many $a \in \mathbb{R}$. Therefore, $K^p(i,j)_a = \Pi'\left(\text{cl}\left(W^p_{\epsilon(i,j),\,a}\right)\right) = \emptyset$ for all but finitely many $a \in \mathbb{R}$. But if $a \notin G^p(i)$ then either $K^p(i)_a = \emptyset$ or $K^p(i)_a$ is a Hausdorff limit of the sets $K^p(i,j)_a$. So $K^p(i)_a = \emptyset$ for all but countably many $a \in \mathbb{R}$. By a similar argument using each $G^p$, we then get that $K^p_a = \emptyset$ for all but countably many $a \in \mathbb{R}$. It follows that the set $G$ is countable. $\square$

**3.12. Corollary.** *Let $X \subseteq \mathbb{R}^m$ be a $T^\infty$-set, and let $1 \leq k \leq m$. Then the set*
$$B := \{\, a \in \mathbb{R}^k :\ \text{cl}(X_a) \neq \text{cl}(X)_a \,\}$$
*has empty interior.*

*Proof.* It suffices to show that the lemma holds with $X \cap \left((-R,R)^k \times \mathbb{R}^{m-k}\right)$ in place of $X$, for each $R > 0$; so we assume that $\Pi_k(X)$ is bounded, where $\Pi_k : \mathbb{R}^m \longrightarrow \mathbb{R}^k$ is the projection onto the first $k$ coordinates. For each $a \in B$ there is a box $U \subseteq \mathbb{R}^{m-k}$ such that $\operatorname{cl}(X_a) \cap U = \emptyset$, but $\operatorname{cl}(X)_a \cap U \neq \emptyset$. Hence $B = \bigcup_U B_U$, where $U$ ranges over all rational boxes in $\mathbb{R}^{m-k}$ and
$$B_U := \{\, a \in \mathbb{R}^k : \operatorname{cl}(X_a) \cap U = \emptyset,\ \operatorname{cl}(X)_a \cap U \neq \emptyset \,\}.$$
Each $B_U$ is contained in the frontier of the bounded $T^\infty$-set $\Pi_k\left(X \cap (\mathbb{R}^k \times U)\right)$. So by the previous proposition $B_U \subseteq Y_U$ for some $T^\infty$-set $Y_U$ with empty interior. Applying Remark 3.10(i) we conclude that $B$ has empty interior. □

**3.13. Proposition.** *If $X \subseteq I^m$ is a $T^\infty$-set, then so is $I^m \setminus X$.*

*Proof.* Let $X \subseteq I^m$ be a $T^\infty$-set. We establish the following two statements by induction on $m$:

(I)$_m$ If $\operatorname{int}(X) = \emptyset$, then $X$ can be partitioned into finitely many $T^\infty$-sets $G_1, \ldots, G_K$ in such a way that for each $k \in \{1, \ldots, K\}$ there is a strictly increasing map $\iota : \{1, \ldots, m-1\} \longrightarrow \{1, \ldots, m\}$ such that $G_k$ is the graph of a continuous function with domain $\Pi_\iota(G_k)$.

(II)$_m$ Each component of $X$ is a $T^\infty$-set, and the complement $I^m \setminus X$ is a $T^\infty$-set.

(The meaning of (I)$_1$ is clear from the conventions on $\iota$ made in the introduction.) The case $m = 1$ follows from Proposition 3.4; so assume $m > 1$ and that the two statements hold for lower values of $m$. First we establish the following

**Claim.** *Assume that there is a $T^\infty$-set $Z \subseteq I^m$ with empty interior such that $X \subseteq Z$ and (I)$_m$ and (II)$_m$ hold with $Z$ in place of $X$. Then (I)$_m$ and (II)$_m$ hold.*

*Proof.* Let $G_1, \ldots, G_K$ be as in (I)$_m$ with $Z$ in place of $X$. Clearly (I)$_m$ then also holds for $X$, since each $G_i \cap X$ is the graph of a continuous function. Since the $G_i$ partition $Z$ and (II)$_m$ holds with $Z$ in place of $X$, it suffices to prove for each $i$ that the components of $G_i \cap X$ and the set $G_i \setminus X$ are $T^\infty$-sets. We may therefore reduce to the case that $K = 1$, that is, there is a strictly increasing map $\iota : \{1, \ldots, m-1\} \longrightarrow \{1, \ldots, m\}$ such that $Z$ is the graph of a continuous function with domain $\Pi_\iota(Z)$. By the inductive hypothesis, each component of $\Pi_\iota(X)$ and the set $\Pi_\iota(Z) \setminus \Pi_\iota(X)$ are $T^\infty$-sets, so the claim follows.

We now return to the proof of the theorem; there are two cases to consider.

**Case 1:** $X$ has empty interior. Let $\Pi$ denote the projection on the first $m-1$ coordinates. Consider the $T^\infty$-sets
$$C_i := \{\, a \in I^{m-1} : |X_a| \geq i \,\} \quad \text{for } i \in \mathbb{N}.$$
By (II)$_{m-1}$ the sets
$$D_i := C_{i+1} \setminus C_i = \{\, a \in I^{m-1} : |X_a| = i \,\}$$
are also $T^\infty$-sets, and by Proposition 3.4 there is an $N \in \mathbb{N}$ such that $C_i = C_{N+1}$ for all $i > N$. The inductive hypothesis therefore allows to reduce to the following two subcases.

**Subcase $C_{N+1} = \emptyset$:** then $|X_a| \leq N$ for every $a \in \Pi(X)$. For $1 \leq j \leq i \leq N$ we define the T$^\infty$-sets
$$X_{i,j} := \{ (a,y) \in C_i \times I : y \text{ is the } j^{\text{th}} \text{ element of } X_a \},$$
$$S_{i,j} := \{ a \in C_i : |(\text{cl}(X_{i,j}))_a| \geq 2 \},$$
and put $S := \bigcup_{1 \leq j \leq i \leq N} S_{i,j}$. (Here we are using the fact that the collection of T$^\infty$-sets is closed under taking topological closure.) (I)$_m$ holds with $X \setminus (S \times I)$ in place of $X$ by construction, and the corresponding (II)$_m$ then follows easily from the inductive hypothesis (and since the order $<$ is definable in $\widetilde{\mathbb{R}}$). Note that
$$S_{i,j} \subseteq \{ a \in \mathbb{R}^{m-1} : \text{cl}((X_{i,j})_a) \neq \text{cl}(X_{i,j})_a \};$$
it follows from Corollary 3.12 that each $S_{i,j}$, and hence $S$, has empty interior. Therefore (I)$_{m-1}$ and (II)$_{m-1}$ hold with $S$ in place of $X$ by the inductive hypothesis, and so (I)$_m$ and (II)$_m$ hold with $S \times I$ in place of $X$. The claim implies now that (I)$_m$ and (II)$_m$ also hold with $X \cap (S \times I)$ in place of $X$.

**Subcase $N = 0$:** then $C_1 = \Pi(X)$. By assumption every fiber $X_a \subseteq I$ with $a \in C_1$ is infinite and hence (by 3.4 again) contains an interval. Since $X$ has empty interior, it follows from Remark 3.10(ii) that $C_1$ has empty interior. (I)$_m$ and (II)$_m$ now follow from the claim by a similar argument as in the previous subcase (with $C_1$ in place of $S$).

**Case 2:** $X$ has nonempty interior. By Lemma 3.11 there is a closed T$^\infty$-set $Y \subseteq I^m$ such that $\text{bd}(X) \subseteq Y$ and $Y$ has empty interior. By Case 1 applied to $Y$, both (I)$_m$ and (II)$_m$ hold with $Y$ in place of $X$. Note that if $C$ is a component of $I^m \setminus Y$ and $C \cap X \neq \emptyset$, then $C \subseteq X$. It follows that each component of $I^m \setminus Y$ is either contained in $X \setminus Y$ or is disjoint from $X \cup Y$. On the other hand, by the claim the statements (I)$_m$ and (II)$_m$ hold with $X \cap Y$ in place of $X$. Thus (II)$_m$ follows easily. □

**3.14. Corollary.** *The collection of T$^\infty$-subsets of $I^m$, $m \in \mathbb{N}$, forms an o-minimal structure on $I$.* □

**3.15. Remark.** Let $(\Lambda_n)_{n \in \mathbb{N}}$ be any system of collections $\Lambda_n$ of subsets of $\mathbb{R}^n$. Replacing throughout the property (for subsets of $\mathbb{R}^n$) of being Pfaffian by the property of belonging to $\Lambda_n$, and adapting the definition of T$^\infty$-set in 3.2 accordingly, all the above from 3.2 onward goes through provided $(\Lambda_n)$ satisfies the following properties:

(i) all subsets of $\mathbb{R}^n$ which are definable in $\widetilde{\mathbb{R}}$ belong to $\Lambda_n$;
(ii) points (1)–(4) of the remark after Definition 3.1 hold with "belongs to $\Lambda_{...}$" in place of "is Pfaffian";
(iii) every $A \in \Lambda_n$ is a closed subset of a definable set $B \subseteq \mathbb{R}^n$;
(iv) for $1 \leq k \leq n$, $A \in \Lambda_n$ and $x \in \mathbb{R}^k$ there is an $N \in \mathbb{N}$ such that each fiber $A_x$ has at most $N$ components;
(v) every $A \in \Lambda_n$ is a finite union of manifolds $B_1, \ldots, B_l \in \Lambda_n$ that are in standard position.

## 4. Pfaffian Closure

Let $\widetilde{\mathbb{R}}$ be an expansion of the real field; in this section "definable" means "definable in $\widetilde{\mathbb{R}}$" unless indicated otherwise. Let $\mathcal{L}(\widetilde{\mathbb{R}})$ be the collection of all Rolle leaves of nonsingular definable 1-forms of class $C^1$ on $\mathbb{R}^n$ (for various $n$). We write $\widetilde{\mathbb{R}}_{\mathcal{L}}$ for the expansion of $\widetilde{\mathbb{R}}$ by all $L \in \mathcal{L}(\widetilde{\mathbb{R}})$.

**4.1. Theorem.** *If $\widetilde{\mathbb{R}}$ is o-minimal, then the structure $\widetilde{\mathbb{R}}_{\mathcal{L}}$ is o-minimal.*

*Proof.* Let $\mathcal{T}_m$ be the collection of all T$^\infty$-sets $X \subseteq I^m$ as defined over $\widetilde{\mathbb{R}}$ in the previous section. By Corollary 3.14, the collection $(\mathcal{T}_m)_m$ forms an o-minimal structure on $I$. Let $\tau_m : \mathbb{R}^m \longrightarrow (0,1)^m$ be the (definable) homeomorphism given by

$$\tau_m(x_1, \ldots, x_m) := \left( \frac{x_1}{1+x_1^2}, \ldots, \frac{x_m}{1+x_m^2} \right),$$

and let $\mathcal{S}_m$ be the collection of sets $\tau_m^{-1}(X)$ with $X \in \mathcal{T}_m$. Then the collection $\mathcal{S} = (\mathcal{S}_m)_m$ gives rise to an o-minimal expansion $\widetilde{\mathbb{R}}_T$ of $\widetilde{\mathbb{R}}$. A routine argument shows that the graphs of addition and multiplication belong to $\mathcal{S}$. But every $L \in \mathcal{L}(\widetilde{\mathbb{R}})$ is definable in $\widetilde{\mathbb{R}}_T$: if $L$ is a Rolle leaf of the definable 1-form $\omega$ of class $C^1$ on $\mathbb{R}^n$ with $S(\omega) = \emptyset$, say, then $\tau_n(L)$ is a Rolle leaf of the pullback $(\tau_n^{-1})^*\omega$. Since $\tau_n(L)$ is a Pfaffian set, it follows from Lemma 3.3 that $\tau_n(L)$ is a T$^\infty$-set and hence definable in $\widetilde{\mathbb{R}}_T$. Therefore $L$ is definable in $\widetilde{\mathbb{R}}_T$, and the theorem is proved. □

**4.2. Lemma.** *If $\widetilde{\mathbb{R}}$ is o-minimal, then every Pfaffian set is definable in $\widetilde{\mathbb{R}}_{\mathcal{L}}$.*

*Proof.* Let $\omega_1, \ldots, \omega_q$ be nonsingular definable 1-forms of class $C^1$ on some open definable set $U \subseteq \mathbb{R}^n$. Let $L_i$ be a Rolle leaf of $\omega_i = 0$ for each $i$, and let $W := A \cap L_1 \cap \cdots \cap L_q$ with $A \subseteq U$ definable. By Lemma 2.1, we may assume that $A$ is a definable $C^1$ cell and that there is a $J \subseteq \{1, \ldots, q\}$ such that $(\omega_i)_{i \in J}$ is transverse to $A$. Lemma 1.6 now gives that every component of $W$ is a component of $A \cap \bigcap_{i \in J} L_i$, so we may even assume that $(\omega_1, \ldots, \omega_q)$ is transverse to $A$. As in the proof of Lemma 2.5 we now reduce further to the case that $U = A = \mathbb{R}^n$. But then $W$ is definable in $\widetilde{\mathbb{R}}_{\mathcal{L}}$. □

In particular, if $\omega$ is a definable 1-form of class $C^1$ on some definable open set $U \subseteq \mathbb{R}^n$, and $L$ is a Rolle leaf of $\omega = 0$, then $L$ is definable in $\widetilde{\mathbb{R}}_{\mathcal{L}}$.

**4.3. Definition.** An expansion $\widetilde{\mathbb{R}}$ of the real field is **Pfaffian closed** if every $L \in \mathcal{L}(\widetilde{\mathbb{R}})$ is definable in $\widetilde{\mathbb{R}}$. Any expansion $\widetilde{\mathbb{R}}$ of the real field admits a smallest expansion $\mathcal{P}(\widetilde{\mathbb{R}})$ which is Pfaffian closed: just repeat the process of adding all Rolle leaves as above to obtain a chain of expansions defined by $\widetilde{\mathbb{R}}_0 := \widetilde{\mathbb{R}}$ and $\widetilde{\mathbb{R}}_{k+1} := (\widetilde{\mathbb{R}}_k)_{\mathcal{L}}$ for $k \geq 0$. Now with $\mathcal{L} := \bigcup_k \mathcal{L}(\widetilde{\mathbb{R}}_k)$, we put

$$\mathcal{P}(\widetilde{\mathbb{R}}) := \left( \widetilde{\mathbb{R}}, (L)_{L \in \mathcal{L}} \right);$$

we call $\mathcal{P}(\widetilde{\mathbb{R}})$ the **Pfaffian closure** of $\widetilde{\mathbb{R}}$.

The theorem of the introduction is now proved by applying Theorem 4.1 and Lemma 4.2 successively to each of the expansions $\widetilde{\mathbb{R}}_k$ above.

Theorem 4.1 and Lemma 4.2 have an interesting consequence: whenever $\widetilde{\mathbb{R}}$ is o-minimal, the collection of $T^\infty$-sets (over $\widetilde{\mathbb{R}}$) generates the same structure on $\mathbb{R}$ as the smaller collection $\mathcal{L}(\widetilde{\mathbb{R}})$ of Rolle leaves. In other words, we have

**4.4. Corollary.** *If $\widetilde{\mathbb{R}}$ is o-minimal, then every $T^\infty$-set is definable in $\widetilde{\mathbb{R}}_{\mathcal{L}}$.*

*Proof.* Let $X \subseteq \mathbb{R}^m$ be a basic $T^\infty$-set defined from $W \subseteq \mathbb{R}^{k+m+l}$, and let $\Pi : \mathbb{R}^{k+m+l} \longrightarrow \mathbb{R}^{k+m}$ be the projection on the first $k+m$ coordinates. Then $X$ is a **limit** of the projection $\Pi(W)$, that is, for each finite set $F \subseteq \mathbb{R}^m$ there is an $a \in \mathbb{R}^n$ such that $X \cap F = \Pi(W)_a \cap F$. By Lemma 4.2 the Pfaffian set $W$ is definable in $\widetilde{\mathbb{R}}_{\mathcal{L}}$, so $\Pi(W)$ is definable in $\widetilde{\mathbb{R}}_{\mathcal{L}}$. Therefore, since $\widetilde{\mathbb{R}}_{\mathcal{L}}$ is o-minimal, a theorem due to Marker and Steinhorn [9] implies that $X$ is definable in $\widetilde{\mathbb{R}}_{\mathcal{L}}$. □

## References

[1] M. Berger and B. Gostiaux, Differential Geometry: Manifolds, Curves, and Surfaces, Springer-Verlag, 1988.
[2] L. van den Dries, o-Minimal structures, in Logic: from Foundations to Applications, W. Hodges et al., eds., Oxford University Press, 1996.
[3] L. van den Dries and C. Miller, Geometric categories and o-minimal structures, Duke Math. J., **84** (1996), 497–540.
[4] M. Karpinski and A. Macintyre, o-Minimal expansions of the real field: A characterization, and an application to pfaffian closure. Preprint, May 1997.
[5] A. Khovanskii, On a class of systems of transcendental equations, Soviet Math. Dokl., **22** (1980), 762–765.
[6] ———, Real analytic varieties with the finiteness property and complex abelian integrals, Funct. Anal. and Appl., **18** (1984), 199–207.
[7] ———, Fewnomials, vol. 88 of Translations of Mathematical Monographs, American Mathematical Society, 1991.
[8] J.-M. Lion and J.-P. Rolin, Feuilletages analytiques réelles et théorème de Wilkie. Preprint, 1996.
[9] D. Marker and C. Steinhorn, Definable types in o-minimal theories, J. Symbolic Logic, **59** (1994), 185–198.
[10] R. Moussu and C. Roche, Théorie de Khovanskii et problème de Dulac, Invent. Math., **105** (1991), 431–441.
[11] ———, Théorèmes de finitude uniformes pour les variétés Pfaffiennes de Rolle, Ann. Inst. Fourier, **42** (1992), 393–420.
[12] A. Pillay and C. Steinhorn, Definable sets in ordered structures. I, Trans. Amer. Math. Soc., **295** (1986), 565–592.
[13] M. Spivak, Differential Geometry, vol. I, Publish or Perish, Inc., Berkeley, 1979.
[14] A. Wilkie, A general theorem of the complement and some new o-minimal structures. Preprint, 1996.

THE FIELDS INSTITUTE, 222 COLLEGE STREET, TORONTO M5T 3J1, ONTARIO, CANADA
*E-mail address*: `pspeisse@fields.utoronto.ca`